\documentclass{article}
\usepackage{amssymb}
\usepackage{amsfonts,amsmath,latexsym}
\usepackage{color}

\catcode`@=11
 \def\.{.\spacefactor\@m}
\catcode`@=12

\begin{document}
\title{Determination of Different Biological Factors on the Base of Dried Blood Spot Technology}
\author{V.~K.~Bozhenko, A.~O.~Ivanov, A.~S.~Mishchenko, \\
A.~A.~Tuzhilin, A.~M.~Shishkin} \maketitle

\section{Introduction}
Determination of different biological factors on the base of dried
blood spot technology has a great practical importance in
investigation of back lands populations, in epidemiological studies
or in special people contingents monitoring, see~\cite{PC},
\cite{PO}. This technology presumes that blood sampling is performed
by patient himself, the sample is spotted on a dry, as a rule,
porous surface (filter paper, cellulose acetate membrane, etc.), and
the posterior transportation to a laboratory, for example, by post
in a standard or a special envelope.

Modern diagnostic equipment gives an opportunity to investigate many
characteristics of dried blood spot, such as metabolites
(see~\cite{Ab}, \cite{Bur}, \cite{PBCC}, \cite{McCann}, and
\cite{Anj}), hormones (see~\cite{Kay}, \cite{Pomelova},
\cite{Levy}), glycated hemoglobin (see~\cite{AnjGee}), and even some
immune system parameters (see~\cite{Shapiro}). The possibility of
DNA and RNA investigations in such samples is of great importance,
since it gives, for example, an opportunity for mass investigations
of socially important infections such as AIDS, hepatitis, etc.

In 1992, in USA a laboratory standard for dried blood spot testing
(DBS) was elaborated, see~\cite{stand}. In 2001, in Russia,
Application Instruction on Alcor Bio Ltd Reagent Kit for
immuno-enzymatic determination of thyrotropic hormone in dried blood
spot of newborn was approved.

Under the dried blood spot technology using, the problem of the
sample spotted volume determination remains  one of the main
practical questions. The existing versions of the technology may
assume spotting of a known blood volume by means of some dosing
device and a posterior elution, or using of special filtering device
for the plasma separation under dried spot preparation,
see~\cite{Ever}.

But there are no any universal method calculating the volume of the
blood spot, which does not use a dosing device. The solution of this
problem gives an opportunity to increase essentially the accuracy of
the results and   to simplify the blood sampling procedure. There is
a series of articles, where the calculations are based on the
concentration of some electrolyte and on a correction of the plasma
volume by the hematocrit value, see~\cite{BBSG}. In the present work
we try to elaborate a universal technology for the spotted volume
calculations.

The solution of this problem is also important for such branches of
Medicine as Catastrophe Medicine and Forensic Laboratory, where
non-standard situations are typical, and, for example, there is no
opportunity to sample patient's blood, and so it is necessary to use
the remains of the patient's blood on some other objects, instead of
the standard blood sample. Such a situation can appear under
investigations of the victims to traffic accidents, or other
catastrophes.

It is well-known that distinct biological indices (analytes) have
distinct variability, see~\cite{Lukich}. We try to use some
mathematical algorithms to pick out a set of blood parameters which
give an opportunity to retrieve the initial volume of the blood
spotted, and use it to calculate exact concentrations of analyts
interesting to a physician. For our analysis we used the database of
biochemical blood parameters obtained in Russian Scientific Center
of Roentgen-Radiology during 1995--2000, which includes more than
$30000$ of patients.

\section{Mathematical Model}
Let us describe the mathematical model of the problem. Let $x_i$,
$1\le i\le m$, stand for the value of the result of the laboratory
analysis on the $i$th molecular compound content. The value $x_i$
obtained as a result of the blood sample analysis depends on two
following parameters at least:  the patient $p$ which is selected
from some collection $\P$, and the volume $\lambda$ of the blood
sample under the analysis. Thus, the value $x_i$ is a function of
two parameters:
 $$
x_i = x_i(p,\lambda).
 $$

The problem is to find a function $f(y_1,y_2,\ldots,y_m)$ whose value
 $$
f\bigl(x_1(p,\lambda),x_2(p,\lambda),\ldots,x_m(p,\lambda)\bigr)
 $$
is close to  $\lambda$ from the statistical point of view.

Notice that due to the uniform distribution of the molecules under
consideration in the blood, a $k$-multiple extension of the volume
must lead to the same increasing of all the indices $x_i$. In other
words, we have $x_i(p,k\,\lambda)=k\,x_i(p,\lambda)$. Therefore, if
$f$ approximates the blood volume, then the following relation must
be valid:
 \begin {multline*}
f\bigl(x_1(p,k\,\lambda),\ldots,x_m(p,k\,\lambda)\bigr)=
f\bigl(k\,x_1(p,\lambda),\ldots,k\,x_m(p,\lambda)\bigr)\approx\\ \approx
k\,f\bigl(x_1(p,\lambda),\ldots,x_m(p,\lambda)\bigr).
 \end {multline*}
This notice is a natural motivation to look for the function  $f$ in
the class of positively homogeneous functions of degree $1$, i.e.,
we assume that the equality
 $$
f(k\,y_1, k\,y_2, \dots, k\,y_m)=k\,f(y_1, y_2, \dots, y_m)
 $$
holds for each positive $k$. Such functions are uniquely defined by
their values at the unit sphere $S^{m-1}$ defined by the equation
$y_1^2+\cdots+y_m^2=1$. By $g$ we denote the restriction of the
function  $f$ onto this unit sphere.

Polynomials form the simplest but rich class of functions. Let us
look for $g$ among the functions which are the restrictions of the
polynomials onto $S^{m-1}$. Our statistical experiments show that it
is enough to consider the polynomials of degree two vanishing at the
origin. In other words, we put $\rho=\sqrt{y_1^2+\cdots+y_m^2}$, and
look for $g$ in the form
 $$
\sum_{i}\alpha_i\frac{y_i}\rho+
\sum_{i\le j}\alpha_{ij}\frac{y_i}\rho\frac{y_j}\rho,
 $$
so the function $f$ is supposed to be in the class
 $$
 f=
\sum_{i}\alpha_i\,y_i+
\sum_{i\le j}\alpha_{ij}y_iy_j/\rho.
 $$

Thus, our problem is to find the coefficients  $\alpha_i$ and
$\alpha_{ij}$ such that the function obtained meets our objectives
as well as possible. To formulate the latter condition
mathematically, let us write down the following objective function.

As we have already mentioned above, the available database gives us
a table of specific values $x_{is}=x_i(p_s,\lambda_s)$. We look for
the function $f$ such that the total squared deviation from the
values $\lambda_s$ is as small as possible. In other words, we have
to find the  $\alpha_i$ and $\alpha_{ij}$ minimizing the objective
function
 \begin {multline*}
h=
\sum_s\bigl(f(x_{1s},\ldots,x_{ms})-\lambda_s\bigr)^2=\\ =
\sum_s\biggl(
\sum_{i}\alpha_i\,x_{is}+
\sum_{i\le j}
\alpha_{ij}x_{is}\frac{x_{js}}{\sqrt{x_{1s}^2+\cdots+x_{ms}^2}}
-\lambda_s\biggr)^2.
 \end {multline*}

Notice that  $h$ considered as a function on  $\alpha_i$ and
$\alpha_{ij}$ is a non-negative quadric. In general position such a
quadric possesses a unique minimum which can be found as a solution
of linear equations system, i.e., from the condition that the
differential of $h$ vanishes.

\section{Application to the specific database}
The above algorithm determining  the volume of a sample for
calculation of the individual values of an arbitrary analyt was
examined on the database of laboratory indices. The best results
were obtained, when we reconstruct the volume by means of the
following analyts: TP, K, Na. The correlation coefficients for the
repaired and true values were 0.95--0.97. The algorithm obtained
gives an opportunity to choose distinct sets of the indices for the
volume reconstruction, that makes the algorithm multipurpose, i.e.
it can be used for analysis of any laboratory blood indices.

The method considered was applied to the specific database in RSCRR.
This database was constructed from $35000$ medical reports
containing biochemical measuring data. We selected the reports
containing the largest number of the biochemical data. So, we
selected the set of $2637$ cases with the next $17$ biochemical data
measured: Chol,    TBil,    DBil, TP,  Alb,  Urea,
     Crea,   ALT, AST,
  Amy,   ALP,  K,
 Ca,   Na,  Fe,
Glu,   LDH.

After calculation of the coefficients $\alpha_i$ and $\alpha_{ij}$
for the function  $f$, we find out the following result: the number
of patients $p_s$ which the inequality
 $$
\frac{|f(x_{1s},\ldots,x_{ms})-\lambda_s|}{\lambda_s}>0,05
 $$
holds for, does not exceed 5\%. This estimate agrees with the
statistical significance of the result.


\end{document}